\documentclass[12pt,a4paper]{amsart}
\usepackage[utf8]{inputenc}	
\usepackage[T1]{fontenc}
\usepackage[in,headings]{fullpage}
\usepackage{amsbsy, amscd, amsfonts, amsmath, amsrefs, amssymb, amsthm}
\usepackage{graphicx}
\usepackage{float}
\usepackage{color}   
\usepackage[colorlinks=true,linkcolor=blue,citecolor=blue,urlcolor=blue]{hyperref}
\def\checked#1{\boxed{#1}}
\def\checked#1{#1}

\newtheorem{theorem}{Theorem}

\newtheorem{proposition}[theorem]{Proposition}
\newtheorem{lemma}[theorem]{Lemma}

\theoremstyle{definition}

\newtheorem{remark}[theorem]{Remark}

\theoremstyle{remark}

\newcommand{\C}{\mathbf{C}}
\newcommand{\Z}{\mathbf{Z}}

\newcommand{\R}{\mathbf{R}}

\renewcommand{\Re}{\mathop{\mathrm{Re}}\nolimits}

\newcommand{\Rzeta}{\mathop{\mathcal R }\nolimits}

\newfont{\cmbsy}{cmbsy10}
\newfont{\cmmib}{cmmib10}
\newcommand{\Orden}{\mathop{\hbox{\cmbsy O}}\nolimits}

\begin{document}

\title{Integral Representations of Riemann auxiliary function.}
\author[Arias de Reyna]{J. Arias de Reyna}
\address{%
Universidad de Sevilla \\ 
Facultad de Matem\'aticas \\ 
c/Tarfia, sn \\ 
41012-Sevilla \\ 
Spain.} 



\email{arias@us.es, ariasdereyna1947@gmail.com}


\begin{abstract}
We prove that the auxiliary function $\Rzeta(s)$ has the integral representation
\[\Rzeta(s)=-\frac{2^s \pi^{s}e^{\pi i s/4}}{\Gamma(s)}\int_0^\infty
y^{s}\frac{1-e^{-\pi y^2+\pi \omega y}}{1-e^{2\pi \omega y}}\,\frac{dy}{y},\qquad \omega=e^{\pi i/4}, \quad\Re s>0,\]
valid for $\sigma>0$. The function in the integrand $\frac{1-e^{-\pi y^2+\pi \omega y}}{1-e^{2\pi \omega y}}$ is entire. Therefore, no residue is added when we move the path of integration. 
\end{abstract}

\maketitle
\section{Introduction}
The auxiliary function of Riemann is defined by the integral
\[\checked{\Rzeta(s)=\int_{0\swarrow1}\frac{x^{-s} e^{\pi i x^2}}{e^{\pi i x}-
e^{-\pi i x}}\,dx.}\]
The position of the zeros of this  function is connected with the zeros of the Riemann zeta function \cite{A166}.

In Section \ref{S:one} we prove the new integral representation
\begin{equation}\label{E:new}
\Rzeta(s)=-\frac{2^s \pi^{s}e^{\pi i s/4}}{\Gamma(s)}\int_0^\infty
y^{s}\frac{1-e^{-\pi y^2+\pi \omega y}}{1-e^{2\pi \omega y}}\,\frac{dy}{y},\qquad \omega=e^{\pi i/4}, \quad\Re s>0,
\end{equation}
Its main interest is that it gives $\Rzeta(s)$ as a Mellin transform of an entire functions. Therefore, changing the path of integration in this integral does not add residues. This gives  new opportunities to bound $\Rzeta(s)$ without the need to bound a zeta sum. Perhaps useful to prove the Lindelöf hypothesis.

In Section \ref{S:two} we give a new proof, starting from \eqref{E:new},  of the representation integral 
\begin{equation}\label{E:anterior}
\Rzeta(s)=\omega e^{\pi i s/4}\sin\frac{\pi
s}{2}\int_0^{+\infty}\frac{y^{-s}e^{-\pi y^2}}{\sin\pi\omega y}\,dy,
\end{equation}
proved in \cite{A166}.

Section \ref{S:tres} gives a form of \eqref{E:new} for $s$ on the critical line:
\begin{equation}
\Rzeta(\tfrac12+it)=(1+ie^{-\pi t})e^{-2i\vartheta(t)}\int_0^\infty \frac{e^{it\log x}}{\sqrt{x}}e^{-\frac{\pi i}{2}(x^2+x)}\frac{\sin(\frac{\pi}{2}(x^2-x))}{\sin(\pi x)}\,dx,
\end{equation}
where the improper integral is convergent.

\section{First expression of \texorpdfstring{$\Rzeta(s)$}{R(s)} }\label{S:one}

We start from an integral representation given by Gabcke \cite{G}. Namely, \begin{equation}\label{IntGabcke}
\checked{\Rzeta(s)=-2^s \pi^{s/2}e^{\pi i s/4}\int_{-\infty}^\infty \frac{e^{-\pi x^2}H_{-s}(x\sqrt{\pi})}{1+e^{-2\pi\omega x}}\,dx,}
\end{equation}
where $\omega=e^{\pi i/4}$ and  $H_\nu(s)$ denotes the Hermite function as defined in the book by Lebedev \cite{L}*{Ch.~10}. The equivalence of representation \eqref{IntGabcke} with the one in Gabcke \cite{G} is shown in Arias de Reyna \cite{A187} where an alternative proof is also given. 

\begin{proposition}\label{P:1}
For $\sigma=\Re s>0$ we have 
\begin{equation}\label{FirstInt}
\checked{
\Rzeta(s)=-\frac{2^s \pi^{s}e^{\pi i s/4}}{\Gamma(s)}\int_0^\infty
y^{s-1}\frac{1-e^{-\pi y^2+\pi \omega y}}{1-e^{2\pi \omega y}}\,dy.}
\end{equation}
\end{proposition}
\begin{proof}
The function $H_\nu(s)$ is an entire function with power series expansion \cite{L}*{eq.(10.4.3)} (except when $\nu$ is a nonnegative integer in which case $H_\nu(z)$ are the usual Hermite polynomials) 
\begin{equation}
\checked{
H_\nu(z)=\frac{1}{2\Gamma(-\nu)}\sum_{n=0}^\infty (-1)^n \Gamma\Bigl(\frac{n-\nu}{2}\Bigr)\frac{(2z)^n}{n!}.}
\end{equation}
Hence we have for $\Re\nu<0$ (with easy justification)
\[\checked{H_\nu(z)=\frac{1}{2\Gamma(-\nu)}\sum_{n=0}^\infty \int_0^\infty y^{\frac{n-\nu}{2}-1}e^{-y}\,dy \frac{(-2z)^n}{n!}=\frac{1}{2\Gamma(-\nu)}\int_0^\infty
y^{-\frac{\nu}{2}}e^{-y-2z\sqrt{y}}\frac{dy}{y}.}\]
Changing variables, $y$ by $\pi y^2$
\[H_{-s}(x\sqrt{\pi})=\frac{\pi^{\frac{s}{2}}}{\Gamma(s)}\int_0^\infty
y^{s}e^{-\pi y^2-2\pi xy}\frac{dy}{y},\qquad \Re s>0.\]

So, \eqref{IntGabcke} implies that
\[\Rzeta(s)=-\frac{2^s \pi^{s}e^{\pi i s/4}}{\Gamma(s)}\int_{-\infty}^\infty \frac{e^{-\pi x^2}}{1+e^{-2\pi\omega x}}\Bigl(\int_0^\infty
y^{s}e^{-\pi y^2-2\pi xy}\frac{dy}{y}\Bigr)\,dx,\qquad \sigma>0.\]
Since 
\[\int_0^\infty
|y^{s}e^{-\pi y^2-2\pi xy}|\frac{dy}{y}\le \checked{\int_0^\infty
y^{\sigma}e^{-\pi y^2-2\pi xy}\frac{dy}{y}=\frac{\Gamma(\sigma)}{\pi^{\sigma/2}}H_{-\sigma}(x\sqrt{\pi}).}\]
By the asymptotic expansion (see Lebedev \cite{L}*{(10.6.6) and (10.6.7)} for $\sigma$ fixed and  $x\in\R$ with $|x|\to+\infty$ we have \[H_{-\sigma}(x\sqrt{\pi})\sim\begin{cases}  (-x\sqrt{\pi})^{-\sigma},& \text{ $x\to+\infty$};\\
-\frac{\sqrt{\pi}e^{-\pi i\sigma}}{\Gamma(\sigma)}e^{\pi x^2}(-x\sqrt{\pi})^{\sigma-1}, &\text{$x\to-\infty$}.\end{cases}\]
So,
\[\frac{\Gamma(\sigma)}{\pi^{\sigma/2}}\int_{-\infty}^\infty 
\Bigl|\frac{e^{-\pi x^2}}{1+e^{-2\pi\omega x}}H_{-\sigma}(x\sqrt{\pi})\Bigr|\,dx<+\infty,\qquad \sigma>0.\]
By Fubini's Theorem we can change the order of integration. 
\[\Rzeta(s)=-\frac{2^s \pi^{s}e^{\pi i s/4}}{\Gamma(s)}\int_0^\infty
y^{s}e^{-\pi y^2}\Bigl(\int_{-\infty}^\infty \frac{e^{-\pi x^2-2\pi xy}}{1+e^{-2\pi\omega x}}\,dx\Bigr)\frac{dy}{y},\qquad \sigma>0.\]
In Siegel's paper about Riemann's nachlass we find 
\[\checked{\int_{0\nwarrow1}\frac{e^{-\pi i u^2+2\pi i uy}}{e^{\pi i u}-e^{-\pi i u}}\,du=
\frac{1}{1-e^{-2\pi i y}}-\frac{e^{\pi i y^2}}{e^{\pi i y}-e^{-\pi i y}}.}\]
Putting $u=\frac12+\omega^3 x$ we obtain 
\[\checked{\int_{0\nwarrow1}\frac{e^{-\pi i u^2+2\pi i uy}}{e^{\pi i u}-e^{-\pi i u}}\,du=
e^{\pi i y}\int_{-\infty}^\infty \frac{e^{-\pi x^2-2\pi \omega xy}}{1+e^{-2\pi \omega x}}\,dx.}\]
It follows that
\[\checked{\int_{-\infty}^\infty \frac{e^{-\pi x^2-2\pi \omega xy}}{1+e^{-2\pi \omega x}}\,dx
=-e^{-\pi i y}\frac{e^{\pi i y^2}-e^{\pi i y}}{e^{\pi i y}-e^{-\pi i y}}.}\]
Hence, for $\sigma>0$ 
\[\Rzeta(s)=-\frac{2^s \pi^{s}e^{\pi i s/4}}{\Gamma(s)}\int_0^\infty
y^{s}e^{-\pi y^2}\Bigl(-e^{-\pi \omega y}\frac{e^{\pi  y^2}-e^{\pi \omega y}}{e^{\pi \omega y}-e^{-\pi \omega y}}\Bigr)\frac{dy}{y},\qquad \sigma>0.\]
That is equivalent to \eqref{FirstInt}.\end{proof}

\begin{remark}
The function $F(z)=\frac{1-e^{-\pi z^2+\pi \omega z}}{1-e^{2\pi \omega z}}$ appearing in \eqref{FirstInt} is entire. The zeros of the denominator are also zeros of the numerator. For $y\to+\infty$ we have $|F(y)|\sim e^{-\pi y\sqrt{2}}$.

The x-ray (see Figure \ref{x-rayF}) shows that the function is relatively small in $|\arg(z)|<\pi/4$ and in the opposite quadrant. We see that there is a line of zeros along  the lines that separate these quadrants on the other, where the function behaves as $e^{-\pi y^2}$.
\begin{figure}
\begin{center}
\includegraphics[width=\hsize]{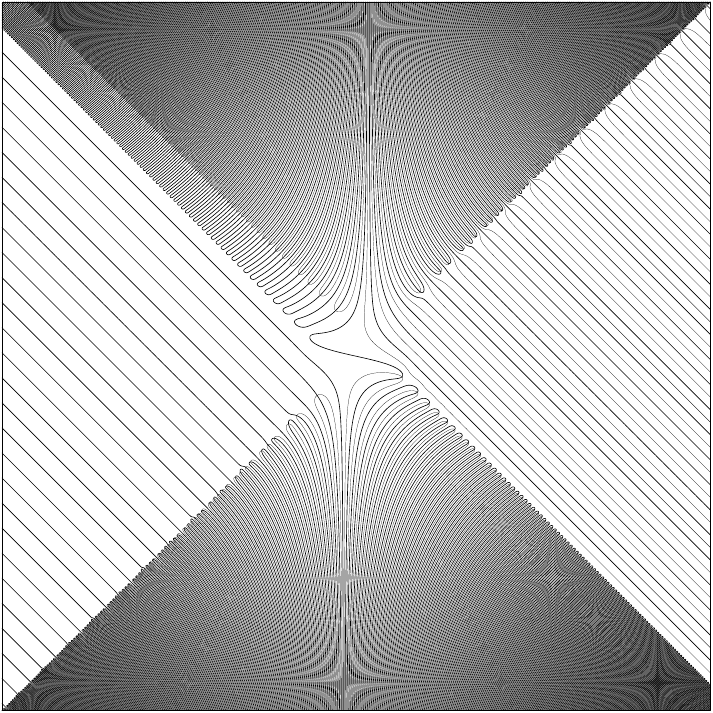}
\caption{x-ray of $\frac{1-e^{-\pi z^2+\pi \omega z}}{1-e^{2\pi \omega z}}$ on $(-10,10)^2$.}
\label{x-rayF}
\end{center}
\end{figure}
\end{remark}

\section{Second integral representation of \texorpdfstring{$\Rzeta(s)$}{R(s)}}\label{S:two}

\begin{proposition}
For $\sigma<0$ we have
\begin{equation}\label{secondInt}
\checked{\Rzeta(s)=-\omega e^{-\pi is/4}(1-e^{\pi i s})\int_0^\infty
y^{-s}\frac{e^{-\pi y^2+\pi \overline\omega y}}{1-e^{2\pi \overline\omega y}}\,dy,\qquad \sigma<0.}
\end{equation}
\end{proposition}
\begin{proof}
For $\sigma>0$, we have the representation \eqref{FirstInt}. For $\sigma>1$ it is easily proved that 
\[\checked{\int_0^\infty y^{s-1}\frac{1}{1-e^{2\pi \omega y}}\,dy=-e^{-\pi i s/4}(2\pi)^{-s}\Gamma(s)\zeta(s).}\]
Therefore, for $\sigma>1$, \eqref{FirstInt} can be written as
\begin{align*}
\Rzeta(s)&=\zeta(s)+\frac{2^s \pi^{s}e^{\pi i s/4}}{\Gamma(s)}\int_0^\infty
y^{s-1}\frac{e^{-\pi y^2+\pi \omega y}}{1-e^{2\pi \omega y}}\,dy\\&=\zeta(s)+\chi(s)e^{-\pi i s/4}(1+e^{\pi i s})\int_0^\infty y^{s-1}\frac{e^{-\pi y^2+\pi \omega y}}{1-e^{2\pi \omega y}}\,dy,\end{align*}
where $\chi(s)$ is the function that appears in the functional equation, it is given by 
\[\chi(s)=\frac{(2\pi)^s}{2\Gamma(s)\cos(\pi s/2)}.\]
In \cite{A166} it is proved that $\zeta(s)=\Rzeta(s)+\chi(s)\Rzeta(1-s)$. Therefore, the above equation is equivalent to 
\[\checked{\overline\Rzeta(1-s)=-e^{-\pi i s/4}(1+e^{\pi i s})\int_0^\infty
y^{s-1}\frac{e^{-\pi y^2+\pi \omega y}}{1-e^{2\pi \omega y}}\,dy}\]
Therefore, putting $1-s$ instead of $s$, we get for $\sigma<0$
\[\checked{\overline\Rzeta(s)=-e^{-\pi i(1-s)/4}(1-e^{-\pi i s})\int_0^\infty
y^{-s}\frac{e^{-\pi y^2+\pi \omega y}}{1-e^{2\pi \omega y}}\,dy}\]
Taking the complex conjugate of both members yields
\[\checked{\Rzeta(\overline s)=-e^{\pi i(1-\overline s)/4}(1-e^{\pi i \overline s})\int_0^\infty y^{-\overline s}\frac{e^{-\pi y^2+\pi \overline\omega y}}{1-e^{2\pi \overline\omega y}}\,dy}.\]
We get \eqref{secondInt} by putting $\overline s$ instead of $s$. 
\end{proof}
\begin{remark}
The equation \eqref{secondInt} is another way to write \cite{A166}*{eq.~(15)}. That is, it is equivalent to saying that for $\sigma<0$ we have
\begin{equation}
\Rzeta(s)=\omega e^{\pi i s/4}\sin\frac{\pi s}{2}\int_0^\infty\frac{y^{-s}e^{-\pi y^2}}{\sin(\pi\omega y)}\,dy.
\end{equation}
\end{remark}

\section{\texorpdfstring{$\Rzeta(s)$}{R(s)} at the critical line}\label{S:tres}

Let $F(z)$ be the function in the integrand of equation \eqref{FirstInt}, that is,
\begin{equation}
F(z):=\frac{1-e^{-\pi z^2+\pi \omega z}}{1-e^{2\pi \omega z}}.
\end{equation}
We need several lemmas on this function.
\begin{lemma}\label{L:231231-2}
There is an absolute constant $C$ such that for $\delta>0$ 
\[|F(\delta+\omega x)|\le C/\delta, \qquad x\ge0.\]
\end{lemma}
\begin{proof}
By definition 
\[F(\delta+\omega x)=\frac{1-e^{-\pi(\delta^2+2\delta\omega x+ix^2)+\pi\omega\delta+\pi i x}}{1-e^{2\pi\delta\omega+2\pi i x}}.\]
And 
\[|1-e^{2\pi\delta\omega+2\pi i x}|\ge e^{\pi\sqrt{2}\delta}-1\ge \pi\sqrt{2}\delta.\]
\[|1-e^{-\pi(\delta^2+2\delta\omega x+ix^2)+\pi\omega\delta+\pi i x}|\le 1+e^{-\pi\delta^2-\pi\sqrt{2}\delta x+\pi\delta/\sqrt{2}}\le 1+e^{-\pi\delta^2+\pi\delta/\sqrt{2}}\le 1+e^{\pi/8}.\]
It follows that 
\[|F(\delta+\omega x)|\le \frac{1+e^{\pi/8}}{\pi\sqrt{2}\delta}\le \frac{2}{3\delta}.\qedhere\]
\end{proof}

\begin{lemma}\label{L:4}
There is an absolute constant $C$ such that for $R>8$ and $0<x<R$ we have $|F(R+ix)|\le C R$.
\end{lemma}
\begin{proof}
In any case we have $|e^{-\pi z^2+\pi\omega z}|\le 1$. To see it, notice that 
\[-\pi z^2+\pi\omega z=-\pi (R^2+2iR x-x^2)+\pi\omega (R+ix).\]
So, 
\[|e^{-\pi z^2+\pi\omega z}|=e^{-\pi(R^2-x^2)+\frac{\pi}{\sqrt{2}}(R-x)}=e^{-\pi(R-x)(R+x+2^{-1/2})}\le 1.\]
Therefore, for $|1-e^{2\pi\omega z}|\ge\frac{1}{R}$ we have $|F(R+ix)|\le2R$. Next, assume that $|1-e^{2\pi\omega z}|<\frac{1}{R}$, and try to prove that, in this case, we also have $|F(R+ix)|\le2R$.

For $|w|<1/4$, we have  
\[\frac12<\Bigl|\frac{e^w-1}{w}\Bigr|<4.\]
If $|1-e^w|<\epsilon<1/8$ we have, for some $n\in\Z$ 
\[|w-2n\pi i|=|\log(1-(1-e^w))|\le\sum_{k=1}^\infty \frac{\epsilon^k}{k}<\frac14.\]
Therefore, the points where $|1-e^w|<\epsilon$, are of the form 
$w=2n\pi i +u$, with $|u|<2\epsilon$ and $n\in\Z$. 
Therefore, for $R>8$, $|1-e^{2\pi\omega z}|<1/R$ implies $2\pi \omega z=2n\pi i +u$, with $|u|<2/R$. We will have  $z=n\omega+\frac{u}{2\pi\omega}$, and $e^{2\pi\omega z}=e^u$.

Then 
\[|1-e^{-\pi z^2+\pi\omega z}|=|1-e^{-\pi (i n^2+\frac{nu}{\pi}+\frac{u^2}{4\pi^2 i})+\pi i n+\frac{u}{2}}|=|1-e^{-nu-\frac{u^2}{4\pi i}+\frac{u}{2}}|\]
Since $z=R+ix=n\omega+\frac{u}{2\pi\omega}$, we have $R=\frac{n}{\sqrt{2}}+\Re\frac{u}{2\pi\omega}$,
then $n=\sqrt{2}R-\sqrt{2}\Re\frac{u}{2\pi\omega}$. It follows that 
\[\Bigl|-nu-\frac{u^2}{4\pi i}+\frac{u}{2}\Bigr|\le \Bigl|-\sqrt{2}u R+u\sqrt{2}\Re\frac{u}{2\pi\omega}+\frac{i u^2}{4\pi}+\frac{u}{2}\Bigr|\le 2\sqrt{2}+\frac{4}{\sqrt{2}\pi R^2}+\frac{1}{\pi R^2}+\frac{1}{R}<4.\]
And there is a constant such that $|(e^z-1)/z|\le C$ for $|z|\le 4$. Therefore, we have for $|1-e^{2\pi\omega z}|<1/R$
\[\Bigl|\frac{1-e^{-\pi z^2+\pi\omega z}}{1-e^{2\pi\omega z}}\Bigr|\le C\frac{|-nu-\frac{u^2}{4\pi i}+\frac{u}{2}|}{|1-e^u|}\le C\frac{|-nu-\frac{u^2}{4\pi i}+\frac{u}{2}|}{|u|/2}\le 2C \Bigl|-n-\frac{u}{4\pi i}+\frac{1}{2}\Bigr|.\]
Therefore,
\[|F(R+ix)|\le 2C\Bigl(\sqrt{2}R+\frac{\sqrt{2}}{\pi R}+\frac{1}{2\pi R}+\frac12\Bigr)\le 4CR.\qedhere \]
\end{proof}

\begin{lemma}\label{L:5}
There is an absolute constant $C$ such that 
\[|F(R+ix)|\le C\min(R,(R-x)^{-1}),\qquad \text{for $R\ge 8$ and $0\le x\le R$.}\]
\end{lemma}
\begin{proof}
The first inequality is true by Lemma \ref{L:4}. For the second inequality, notice that 
$R+ix=(R-x)+x(1+i)=\delta+i\omega y$ with $\delta=R-x$ and $y=\sqrt{2}x$. Therefore, by Lemma \ref{L:231231-2}, we have $|F(R+ix)|\le C/(R-x)$. 
\end{proof}

\begin{proposition}\label{P:ocho}
For $s=\frac12+it$ with $t>0$ we have
\begin{equation}\label{Oh!}
\Rzeta(\tfrac12+it)=e^{\frac{\pi i}{4}}\frac{(2\pi)^{\frac12+it}e^{-\frac{\pi t}{2}}}{\Gamma(\frac12+it)}\int_0^\infty \frac{e^{it\log x}}{\sqrt{x}}e^{-\frac{\pi i}{2}(x^2+x)}\frac{\sin(\frac{\pi}{2}(x^2-x))}{\sin(\pi x)}\,dx.
\end{equation}
\end{proposition}
\begin{proof}
Let $s=\frac12+it$, with $t>0$. By Cauchy's Theorem, for $R>0$, we have
\[\int_0^{R+iR}z^{s-1}F(z)\,dz=\int_0^R z^{s-1}F(z)\,dz+\int_R^{R+iR}z^{s-1}F(z)\,dz.\]
For the first of these two integrals, by Proposition \ref{P:1}, we have 
\[\lim_{R\to+\infty}\int_0^R z^{s-1}F(z)\,dz=\int_0^\infty z^{s-1}F(z)\,dz=-\frac{\Gamma(s)}{(2\pi)^s e^{\pi i s/4}}\Rzeta(s).\]
For the second integral, we have the bound 
\[|I_2|:=\Bigl|\int_R^{R+iR}z^{s-1}F(z)\,dz\Bigr|\le \int_0^R|(R+ix)^{-\frac12+it}F(R+ix)|\,dx.\]
We have 
\begin{align*}
|(R+ix)^{-\frac12+it}|&=|\exp((-\tfrac12+it)(\tfrac12\log(R^2+x^2)+i\arctan(x/R))|\le 
(R^2+x^2)^{-1/4} \\ &\le R^{-1/2}.\end{align*}
By Lemma \ref{L:5} we have
\[|I_2|\le C R^{-1/2}\Bigl(\int_0^{R-1/R}\frac{dx}{R-x}+\int_{R-1/R}^R CR\,dx\Bigr)\le C R^{-1/2}(2\log R+C).\]
Therefore, the improper integral converges and 
\[\int_0^\infty (\omega x)^{s-1}F(\omega x)\omega\,dx=-\frac{\Gamma(s)}{(2\pi)^s e^{\pi i s/4}}\Rzeta(s).\]
Notice that 
\[(\omega x)^{s-1}=\exp((-\tfrac12+it)(\log x+\tfrac{\pi i}{4})=\frac{e^{-\pi t/4}}{\sqrt{x}} e^{it\log x-\pi i/8},\]
and
\[F(\omega x)=\frac{1-e^{-\pi i x^2+\pi i x}}{1-e^{2\pi i x}}=-\frac{e^{-\frac{\pi i}{2}(x^2-x)}}{e^{\pi i x}}\frac{\sin(\frac{\pi}{2}(x^2-x))}{\sin(\pi x)}=-e^{-\frac{\pi i}{2}(x^2+x)}\frac{\sin(\frac{\pi}{2}(x^2-x))}{\sin(\pi x)}.\]
Substituting and reordering, we get \eqref{Oh!}. 
\end{proof}

\begin{proposition}
For $t>0$ we have 
\begin{equation}
\Rzeta(\tfrac12+it)=(1+ie^{-\pi t})e^{-2i\vartheta(t)}\int_0^\infty \frac{e^{it\log x}}{\sqrt{x}}e^{-\frac{\pi i}{2}(x^2+x)}\frac{\sin(\frac{\pi}{2}(x^2-x))}{\sin(\pi x)}\,dx
\end{equation}
\end{proposition}
\begin{proof}
By \eqref{Oh!} we only need to show that 
\[\checked{e^{\frac{\pi i}{4}}\frac{(2\pi)^{\frac12+it}e^{-\frac{\pi t}{2}}}{\Gamma(\frac12+it)}=(1+ie^{-\pi t})e^{-2i\vartheta(t)}.}\]
Putting $s=\frac12+it$, we have $t=i(\frac12-s)$. The left-hand side is transformed into 
\[e^{\frac{\pi i}{4}}\frac{(2\pi)^{\frac12+it}e^{-\frac{\pi t}{2}}}{\Gamma(\frac12+it)}=\frac{(2\pi)^s e^{\pi i s/2}}{\Gamma(s)}=(1+e^{\pi i s})\frac{(2\pi)^s}{2\Gamma(s)\cos\frac{\pi s}{2}}.\]
From Titchmarsh \cite{T}*{4.17} we get  $e^{-2i\vartheta(t)}=\chi(s)$,  and $\chi(s)=\frac{(2\pi)^s}{2\Gamma(s)\cos\frac{\pi s}{2}}$.
\end{proof}

\section{Trying to bound \texorpdfstring{$\Rzeta(s)$}{R(s)}}

We have not been able to take advantage of integrals to bound $\Rzeta(s)$. In this section, we expose one of these failed attempts.

\begin{lemma}
For  $s=\frac12+it$ with $t\in \R$ and $a\in\R$ fixed,   we have
\[\lim_{R\to+\infty}\int_R^{R+ia}z^{s-1}\frac{1-e^{-\pi z^2+\pi \omega z}}{1-e^{2\pi \omega z}}\,dz=0.\]
\end{lemma}
\begin{proof}
For $|a|<R$ and $z=R+ix$ with $|x|<|a|$ we have 
\begin{multline*}|z^{s-1}|=|\exp((-\tfrac12+it)(\log(R^2+x^2)^{1/2}+i\arctan\tfrac{x}{R}))|\\=
\exp(-\tfrac12\log(R^2+x^2)^{1/2}-t\arctan\tfrac{x}{R})\le R^{-1/2}e^{\frac{\pi |t|}{2}}.\end{multline*}
There is some $R_0(a)$ such that for  $R>R_0(a)$ 
\begin{multline*}|1-e^{-\pi z^2+\pi\omega z}|=|1-e^{-\pi(R^2+2iRx -x^2)+\pi\frac{1+i}{\sqrt{2}}(R+ix)}|\\
\le 1+e^{-\pi(R^2-x^2)+\frac{\pi R}{\sqrt{2}}-\frac{\pi x}{\sqrt{2}}}
\le 1+e^{\pi a^2+\frac{\pi |a|}{\sqrt{2}}}e^{-\pi R^2+\frac{\pi R}{\sqrt{2}}}
\le 2,\end{multline*}
and
\[|1-e^{2\pi\omega z}|=|1-e^{\pi\sqrt{2}(1+i)(R+ix)}|\ge e^{\pi R\sqrt{2}-\pi x\sqrt{2}}-1\ge e^{-\pi \sqrt{2}|a|}e^{\pi R\sqrt{2}}-1\ge e^{\pi R}.\]
Therefore,
\[\Bigl|\int_R^{R+ia}z^{s-1}\frac{1-e^{-\pi z^2+\pi \omega z}}{1-e^{2\pi \omega z}}\,dz\Bigr|\le 2R^{-1/2}e^{\frac{\pi |t|}{2}} e^{-\pi R}|a|.\qedhere\]
\end{proof}
From the Lemma and Cauchy's Theorem it follows directly the next Proposition.
\begin{proposition}\label{P:2}
For any $a\in\C$ we have 
\begin{equation}\label{FirstInt2}
\Rzeta(s)=-\frac{2^s \pi^{s}e^{\pi i s/4}}{\Gamma(s)}\int_{\Gamma_a}
z^{s-1}\frac{1-e^{-\pi z^2+\pi \omega z}}{1-e^{2\pi \omega z}}\,dz,
\end{equation}
where $\Gamma_a$ is the path composed of the segment $[0,a]$ and the half-line
$[a,a+\infty)$ parallel to the real axis.
\end{proposition}

\begin{proposition}
For any $t\in\R$ and $\rho>0$ we have
\begin{equation}
\checked{\Rzeta(\tfrac12+it)=\frac{e^{-\pi t/2}}{\Gamma(\frac12+it)}e^{it\log t}(K_1+K_2),}
\end{equation}
where
\begin{equation}\label{EK1}
\begin{aligned}
K_1&=\omega t^{1/2}\int_0^{\rho}
\frac{e^{it\log x-itx}}{\sqrt{x}}\frac{1-e^{-\frac{it^2x^2}{4\pi}+\frac{itx}{2}}}{1-e^{-itx}}\,dx\\&=\omega e^{-it\log t}\int_0^{\rho t}\frac{e^{it\log x-ix}}{\sqrt{x}}\frac{1-e^{-\frac{ix^2}{4\pi}+\frac{ix}{2}}}{1-e^{-ix}}\,dx,
\end{aligned}
\end{equation}
\begin{equation}\label{EK2}
K_2=(\omega\rho t)^{1/2}e^{it(\log\rho-\rho)}\int_0^\infty\frac{e^{\frac\pi4 t+it\log(\omega+x)-t\rho\omega x}}{\sqrt{\omega+x}}\frac{1-e^{-\frac{\rho^2t^2}{4\pi}(\omega+x)^2+\frac{\rho t}{2}(i+\omega x)}}{1-e^{-\rho t(i+\omega x)}}\,dx.
\end{equation}
\end{proposition}
\begin{proof}
Taking $a=\rho t\omega/2\pi$ in Proposition \ref{P:1} we get 
\[\Rzeta(\tfrac12+it)=-\frac{2^s \pi^{s}e^{\pi i s/4}}{\Gamma(s)}\Bigl(\int_0^{\rho t\omega/2\pi}
z^{s-1}\frac{1-e^{-\pi z^2+\pi \omega z}}{1-e^{2\pi \omega z}}\,dz+
\int_{\rho t\omega/2\pi}^{\rho t\omega/2\pi+\infty}
z^{s-1}\frac{1-e^{-\pi z^2+\pi \omega z}}{1-e^{2\pi \omega z}}\,dz\Bigr)\]
In the first integral, putting $z=\frac{t\omega}{2\pi} x $ with $0<x<\rho$ we get
\begin{multline*}\int_0^{\rho t\omega/2\pi}
z^{s-1}\frac{1-e^{-\pi z^2+\pi \omega z}}{1-e^{2\pi \omega z}}\,dz=\int_0^{\rho}
\frac{(\frac{t}{2\pi})^{-1/2}e^{-\frac{\pi t}{4}}e^{it\log\frac{tx}{2\pi}-\frac{\pi i}{8}}}{\sqrt{x}}\frac{1-e^{-\frac{it^2x^2}{4\pi}+\frac{itx}{2}}}{1-e^{itx}}\frac{t\omega}{2\pi}\,dx\\
=-\Bigl(\frac{t}{2\pi}\Bigr)^{1/2}e^{-\frac{\pi t}{4}}e^{it\log \frac{t}{2\pi}+\frac{\pi i}{8}}
\int_0^{\rho}
\frac{e^{it\log x-itx}}{\sqrt{x}}\frac{1-e^{-\frac{it^2x^2}{4\pi}+\frac{itx}{2}}}{1-e^{-itx}}\,dx\end{multline*}
In the second integral taking $z=\frac{\rho t}{2\pi}(\omega+x)$ with $0<x<+\infty$, we get
\begin{multline*}
\int_{\frac{\rho t\omega}{2\pi}}^{\frac{\rho t\omega}{2\pi}+\infty}\mskip-4mu
z^{s-1}\frac{1-e^{-\pi z^2+\pi \omega z}}{1-e^{2\pi \omega z}}\,dz=
\Bigl(\frac{\rho t}{2\pi}\Bigr)^{\frac12}e^{it\log\frac{\rho t}{2\pi}}\mskip-10mu\int_0^\infty\frac{e^{it\log(\omega+x)}}{\sqrt{\omega+x}}\frac{1-e^{-\frac{\rho^2t^2}{4\pi}(\omega+x)^2+\frac{\rho t}{2}(i+\omega x)}}{1-e^{\rho t(i+\omega x)}}\,dx\\
=-\Bigl(\frac{\rho t}{2\pi}\Bigr)^{\frac12}e^{it\log\frac{\rho t}{2\pi}-i\rho t}\int_0^\infty\frac{e^{it\log(\omega+x)-t\rho\omega x}}{\sqrt{\omega+x}}\frac{1-e^{-\frac{\rho^2t^2}{4\pi}(\omega+x)^2+\frac{\rho t}{2}(i+\omega x)}}{1-e^{-\rho t(i+\omega x)}}\,dx.
\end{multline*}
Joining this, we get the equality
\[\Rzeta(\tfrac12+it)=\frac{e^{-\pi t/2}}{\Gamma(\frac12+it)}e^{it\log t}(K_1+K_2),\]
where $K_1$ and $K_2$ are given by \eqref{EK1} and \eqref{EK2}.
\end{proof}
\begin{remark}
By Stirling expansion, we have for $t\to+\infty$
\[\frac{e^{-\pi t/2}}{\Gamma(\frac12+it)}=\frac{1}{\sqrt{2\pi}}e^{-(it\log t-it)}(1+\Orden(t^{-1})).\]
Therefore, the Lindelöf hypothesis is equivalent to $K_1+K_2\ll t^\varepsilon$. We will prove that $K_2=\Orden(1)$, but this does not seem to bring us any closer to Lindelof's hypothesis.\end{remark}

\begin{proposition}
For $t>0$, $\rho>1$ and $\rho t>\pi$ we have the bound 
\begin{equation}
|K_2|\le\frac{2^{5/4}}{|\sin(\rho t/2)|}\Bigl(\frac{\pi\rho}{\rho-1}\Bigr)^{1/2}.
\end{equation}
\end{proposition}
\begin{proof}
We have
\[|1-e^{-\rho t(i+\omega x)}|\ge 1-e^{-\rho t x/\sqrt{2}},\]
and
\begin{multline*}
|1-e^{-\rho t(i+\omega x)}|=|e^{\rho i t/2}-e^{-\rho i t/2-\rho t\omega x}|=
|e^{\rho i t/2}-e^{-\rho i t/2}+e^{-\rho i t/2}(1-e^{-\rho t\omega x})|\\
2|\sin(\rho t/2)|-(1-e^{-\rho t x/\sqrt{2}}).
\end{multline*}
From both we derive that in general 
\[|1-e^{-\rho t(i+\omega x)}|\ge |\sin(\rho t/2)|.\]
Also,
\[|1-e^{-\frac{\rho^2t^2}{4\pi}(\omega+x)^2+\frac{\rho t}{2}(i+\omega x)}|\le 
1+e^{-\frac{\rho^2t^2}{4\pi}\sqrt{2}x-\frac{\rho^2t^2}{4\pi}x^2+\frac{\rho t }{2\sqrt{2}}x}.\]
The exponent is maximum at $x=\frac{\pi-\rho t}{\sqrt{2}\,\rho t}$ Assuming $\rho t>\pi$ we have  $-\frac{\rho^2t^2}{4\pi}\sqrt{2}x-\frac{\rho^2t^2}{4\pi}x^2+\frac{\rho t }{2\sqrt{2}}x<0$ for $x>0$ and therefore
\[|1-e^{-\frac{\rho^2t^2}{4\pi}(\omega+x)^2+\frac{\rho t}{2}(i+\omega x)}|\le
2.\]

It follows that 
\[|K_2|\le (\rho t)^{1/2}\frac{2}{|\sin(\rho t/2)|}
\int_0^\infty\frac{e^{\frac{\pi t}{4}-t\arctan\frac{1}{1+x\sqrt{2}}-\frac{t\rho x}{\sqrt{2}}}}{\sqrt{|\omega+x|}}\,dx\]
It is easy to show that for $x>0$ 
\[\frac{\pi}{4}-\arctan\frac{1}{1+x\sqrt{2}}-\frac{\rho x}{\sqrt{2}}\le \frac{1-\rho}{\sqrt{2}}x.\]
Hence assuming $\rho>1$ and $t>0$ we get 
\[\int_0^\infty\frac{e^{\frac{\pi t}{4}-t\arctan\frac{1}{1+x\sqrt{2}}-\frac{t\rho x}{\sqrt{2}}}}{\sqrt{|\omega+x|}}\,dx\le \int_0^\infty\frac{e^{\frac{t(1-\rho) x}{\sqrt{2}}}}{\sqrt{x}}\,dx=\frac{2^{1/4}\sqrt{\pi}}{\sqrt{t(\rho-1)}}.\qedhere\]
\end{proof}

\end{document}